\documentclass[
final
]{dmtcs-episciences}

\usepackage[utf8]{inputenc}
\usepackage{subfigure}

\usepackage{amsmath}

\author[Kenta Noguchi and Carol T. Zamfirescu]{Kenta Noguchi\affiliationmark{1}\thanks{Partially supported by JSPS KAKENHI Grant Number JP21K13831.}
  \and Carol T. Zamfirescu\affiliationmark{2,3}
}
\title[Spanning trees for many different numbers of leaves]{Spanning trees for many different numbers of leaves}
\affiliation{
  Tokyo University of Science, Tokyo, Japan\\
  Ghent University, Ghent, Belgium\\
  Babe\c{s}-Bolyai University, Cluj-Napoca, Roumania}
\keywords{spanning tree, leaf counting, triangulation, polynomial time}
\begin{document}
\publicationdata
{vol. 26:3}
{2024}
{17}
{10.46298/dmtcs.13116}
{2024-02-26; 2024-02-26; 2024-10-31}
{2024-10-31}
\maketitle
\begin{abstract}
Let $G$ be a connected graph and $L(G)$ the set of all integers $k$ such that $G$ contains a spanning tree with exactly $k$ leaves. We show that for a connected graph $G$, the set $L(G)$ is contiguous. It follows from work of Chen, Ren, and Shan that every connected and locally connected $n$-vertex graph -- this includes triangulations -- has a spanning tree with at least $n/2 + 1$ leaves, so by a classic theorem of Whitney and our result, in any plane $4$-connected $n$-vertex triangulation one can find for any integer $k$ which is at least~2 and at most $n/2 + 1$ a spanning tree with exactly $k$ leaves (and each of these trees can be constructed in polynomial time). We also prove that there exist infinitely many $n$ such that there is a plane $4$-connected $n$-vertex triangulation containing a spanning tree with $2n/3$ leaves, but no spanning tree with more than $2n/3$ leaves. 
\end{abstract}

\section{Introduction}

At the 2023 Montreal Graph Theory Workshop, Kenta Ozeki raised several questions regarding the number of leaves in spanning trees of planar 4-connected graphs, and in particular in plane 4-connected triangulations. One such problem was to determine, for a given plane 4-connected triangulation $G$, all integers $k$ for which there exists a spanning tree of $G$ with exactly $k$ leaves. In this note we give a partial answer to this question. 

Let $G$ be a connected graph and $L(G)$ the set of all integers $k$ such that $G$ contains a spanning tree with exactly $k$ leaves. We first prove that in any connected graph $G$, the set $L(G)$ is contiguous. We were surprised to not be able to find this result in the literature. On the one hand, the proof is not difficult, but on the other hand, we believe this result is of independent interest.

By this result, the problem of determining $L(G)$ for a given graph $G$ is equivalent to determining the {\it minimum leaf number} and {\it maximum leaf number} of $G$, i.e.\ the number of leaves in a spanning tree of $G$ with the fewest (resp.\ the most) leaves among all spanning trees of $G$.

The decision problems associated to both aforementioned problems are NP-hard for arbitrary graphs $G$. Thus, determining $L(G)$ is NP-hard. Solis-Oba gave a 2-approximation algorithm for the maximum leaf number problem~\cite{So98}. It was proven by Lu and Ravi~\cite{LR96} that, unless P = NP, there is no constant factor approximation for minimising the number of leaves of a spanning tree. From an optimisation standpoint, minimising the number of leaves is equivalent to maximising the number of non-leaves; perhaps counter-intuitively, for the latter problem there does exist a (linear-time) 2-approximation algorithm based on depth first search, as shown by Salamon and Wiener~\cite{SW08}. 

We now focus on Ozeki's problem asking for $L(G) =: L$ when $G$ is a plane 4-connected triangulation. Since planar 4-connected graphs are hamiltonian~\cite{Tu56}, $\min L = 2$. It remains to determine $\max L$ which seems to be a challenging problem. Let $f(n)$ be the largest integer~$\ell$ such that every plane 4-connected triangulation on $n$ vertices contains a spanning tree with at least $\ell$ leaves.
 
In a given graph $G$, a {\it homeomorphically irreducible spanning tree (HIST)} of $G$ is a spanning tree of $G$ with no vertices of degree~2. 
Chen, Ren, and Shan~\cite{CRS12} showed the following results which will be useful to us. 
Here, a graph $G$ is \emph{locally connected} if for every vertex $v \in V(G)$, the subgraph induced by the neighbourhood $N(v)$ is connected. 

\bigskip

\noindent \textbf{Theorem A} (Theorem 1.1.\ from \cite{CRS12})\textbf{.} \emph{Every connected and locally connected graph with order at least four contains a HIST.}

\bigskip

\noindent \textbf{Corollary B} (Corollary 1.2.\ from \cite{CRS12})\textbf{.} \emph{Let $\Pi$ be a surface (orientable or non-orientable). Then every triangulation of $\Pi$ with at least four vertices contains a HIST.}

\bigskip

It easily follows from Corollary B that every $n$-vertex triangulation of any surface (not necessarily 4-connected) has a spanning tree with at least $\frac{n}{2} + 1$ leaves for all $n\ge 4$. 
So $f(n) \ge \frac{n}{2} + 1$ for $n \ge 4$; we note that this bound can also be inferred from an earlier result of Albertson et al. \cite{ABHT90}. 
Next to our (aforementioned) main theorem, we prove that $f(n) \le \frac{2n}{3}$ for infinitely many $n$. 

Regarding the problem of determining bounds for $f(n)$, on the one hand, there are some results concerning minimum degree conditions. 
Let $g_k(n)$ be the largest integer $\ell$ such that every $n$-vertex graph with minimum degree $k$ contains a spanning tree with at least $\ell$ leaves. 
Linial~\cite{GKS89} conjectured $g_k(n) = \frac{(k-2)n}{(k+1)}+c_k$ for some constant $c_k$. For $k=3$, the tight value $g_3(n) = \frac{n}{4}+2$ is shown in \cite{KW91}, while for $k=4$, the tight value $g_4(n) = \frac{2n}{5}+\frac{8}{5}$ is shown in \cite{GW92, KW91}. (See also the extension of Kleitman and West's result~\cite{KW91} by Karpov~\cite{Ka12}; $g_4(n) = \frac{2n}{5}+2$ holds except for three graphs.) For $k=5$, the tight value $g_5(n) = \frac{n}{2}+2$ is shown in \cite{GW92}. 
For $k \ge 6$, the problem is still open. On the other hand, as far as we know, 
there are few results (or even problems) regarding global or topological conditions, e.g.~connectedness. 

We briefly discuss the link between our problem and domination. 
For a spanning tree of a connected graph $G$, its stem vertices (i.e.\ non-leaves) correspond to a \emph{connected dominating set} in $G$, and our problem is to minimise the size of a connected dominating set, see for instance~\cite{CWY00}. 
For a triangulation $G$ on any surface, the domination number $\gamma(G)$ is at most $\frac{n}{3}$, see~\cite{BHT12} (cf.~\cite{FM15}) and it is conjectured that the tight lower bound is $\frac{n}{4}$ for sufficiently large $n$, see \cite{MT96,PZ09}. 
This gap suggests that it is also difficult to improve the upper bound of the size of a connected domination set. 
(However, the upper bound of the domination number for \emph{plane} triangulations has improved; see \cite{CRR,PYZ16,PYZ20,Sp20}.) 
See also Remark 2 in the last paragraph in Section \ref{sec2}.

\section{Results}
\label{sec2}

For a given tree $T$, let $\ell(T)$ denote the number of leaves in $T$. In our proof, we transform one spanning tree into another via edge exchanges. This technique has been used in the so-called \textit{reconfiguration problem}; see for example~\cite{BI20,ID11} for related problems.

\bigskip

\noindent \textbf{Theorem.} \emph{For any connected graph $G$, the set $L(G)$ is contiguous.}

\bigskip

\begin{proof}
The statement is clearly true for trees, so henceforth we assume $G$ to not be a tree. Let $T$ and $T'$ be arbitrary but fixed spanning trees of $G$. 
Let $E(T') \setminus E(T) = \{ e_0, \ldots, e_{p - 1} \}$ with $p$ positive. 
Put $T_0 := T$ and $T_{i+1} := T_i + e_i - e'_i$ for $i \in \{ 0, \ldots, p - 1 \}$, where $e'_i \in E(C^i) \setminus E(T')$, where $C^i$ is the (unique) cycle in $T_i + e_i$. We need the following two claims.

\bigskip

\noindent \textbf{Claim 1.} $|\ell(T_i) - \ell(T_{i+1})| =: \ell_i \le 2$.

\smallskip

Comparing the degrees of vertices in $T_i$ and $T_{i+1}$, changes can only occur in the end-vertices of $e_{i} = u_{i}v_{i}$ and $e'_i = u'_i v'_i$.
We have
\begin{align*}
d_{T_{i+1}}(u_{i}) &\in \{ d_{T_i}(u_{i}), d_{T_i}(u_{i}) + 1 \},\\
d_{T_{i+1}}(v_{i}) &\in \{ d_{T_i}(v_{i}), d_{T_i}(v_{i}) + 1 \},\\
d_{T_{i+1}}(u'_i) &\in \{ d_{T_i}(u'_i), d_{T_i}(u'_i) - 1 \}, \textrm{and}\\
d_{T_{i+1}}(v'_i) &\in \{ d_{T_i}(v'_i), d_{T_i}(v'_i) - 1 \}.
\end{align*}
Thus, letting $\varepsilon_1, \varepsilon_2 \in \{ -1, 0 \}$ and $\varepsilon_3, \varepsilon_4 \in \{ 0, 1 \}$, by the above argument we have
$$\ell(T_{i+1}) = \ell(T_i) + \sum_{i=1}^4 \varepsilon_i = \ell(T_i) + \varepsilon$$
for $\varepsilon \in \{ -2, - 1, 0, 1, 2 \}.$

\smallskip
Note that if $e_i$ and $e'_i$ share a common vertex, say $v_i = v'_i$, we can make a stronger statement:
$$\ell(T_{i+1}) = \ell(T_i) + \varepsilon_1 + \varepsilon_3 = \ell(T_i) + \varepsilon'$$
for $\varepsilon' \in \{- 1, 0, 1 \}.$ Thus, Claim 1 is proven.
\qed

\bigskip

\noindent \textbf{Claim 2.} \emph{If, in Claim~1, we have $\ell_i = 2$, then there exists a spanning tree $S_i$ of $G$ such that $\ell(T_i) < \ell(S_i) < \ell(T_{i+1})$ when $\ell(T_i) < \ell(T_{i+1})$, and $\ell(T_{i+1}) < \ell(S_i) < \ell(T_i)$ when $\ell(T_i) > \ell(T_{i+1})$.}

\smallskip

In $C^i$, i.e.\ the cycle in $T_i + e_i$, consider a sequence of edges $e_i = e_i^0, e_i^1, \ldots, e_i^{t-1}, e_i^t = e'_i$ (recall that $e'_i \in E(C^i) \setminus E(T')$) such that consecutive edges share a vertex.
Put $T_{i_0} := T_i$ and $T_{i_x} := T_{i_0} + e_i^0 - e_i^x$ for $x \in \{ 1, \ldots, t \}$.
(Here notice that $T_{i_t} = T_{i+1}$.)
Then $|\ell(T_{i_{x-1}}) - \ell(T_{i_x})| \le 1$, which follows from the last remark in the proof of Claim~1. This completes the proof of Claim~2.
\qed

\bigskip

We can now finish the proof of the theorem. Let $T$ and $T'$ be spanning trees of $G$ with the minimum number and the maximum number of leaves, respectively. By Claims~1 and 2, there exists a sequence of spanning trees $T = R_0, R_1, \ldots, R_{q-1}, R_q = T'$ in $G$ such that $|\ell(R_i) - \ell(R_{i+1})| \le 1$ for all $i \in \{ 0, \ldots, q - 1 \}$.
This implies the theorem's statement.
\end{proof}

\bigskip

One might wonder, when $\ell(T) < \ell(T')$, whether there exists a ``monotone'' sequence of spanning trees $T = T_0, T_1, \ldots, T_k = T'$, that is, $|E(T_i) \backslash E(T_{i+1})| = 1$ and $\ell(T_i) \le \ell(T_{i+1})$ for all $i \in \{ 0, \ldots, k-1 \}$. In general, there are counterexamples, e.g.\ the following, where we define a graph $G$ and two spanning trees $T$ and $T'$ in $G$: 
\begin{align*}
V(G) &= \{ v_1, \ldots, v_8\},\\
E(G) &= \{ v_1v_2, v_1v_3, v_1v_4, v_1v_5, v_1v_6, v_2v_7, v_3v_7, v_4v_8, v_5v_8, v_6v_7, v_6v_8 \},\\
E(T) &= \{ v_1v_6, v_2v_7, v_3v_7, v_4v_8, v_5v_8, v_6v_7, v_6v_8 \}, \textrm{and}\\
E(T') &= \{ v_1v_2, v_1v_3, v_1v_4, v_1v_5, v_1v_6, v_6v_7, v_6v_8 \}.
\end{align*}
Now $T'$ attains the maximum leaf number of $G$ and $\ell(T) = 5 < 6 = \ell(T')$ but every edge exchange applied to $T$ decreases the number of leaves. However, we could not answer the following question.

\bigskip

\noindent \textbf{Problem.} \textit{Consider a connected graph $G$. Let $T$ and $T'$ be spanning trees of $G$ with the minimum number and the maximum number of leaves, respectively. Is it then true that there exists a sequence of spanning trees $T = T_0, T_1, \ldots, T_k = T'$ such that $|E(T_i) \backslash E(T_{i+1})| = 1$ and $\ell(T_i) \le \ell(T_{i+1})$ for all $i \in \{ 0, \ldots, k-1 \}$?}

\bigskip

We now mention some consequences of our theorem. Recall that planar 4-connected graphs are hamiltonian~\cite{Tu56}, so they contain a hamiltonian path, i.e.\ a spanning tree with exactly two leaves (we call graphs containing a hamiltonian path \textit{traceable}).

\bigskip

\noindent \textbf{Corollary.} \emph{Let $G$ be a traceable graph. Then $G$ contains a spanning tree with exactly $k$ leaves for all $k \in \{ 2, \ldots, \max L(G) \}$. In particular, this holds for planar $4$-connected graphs.}

\bigskip

We also remark that if one can find in a connected and locally connected graph $G$ on $n$~vertices a spanning tree with $p < \frac{n}{2}+1$ leaves, then one can find a spanning tree of $G$ with $r$ leaves for every integer $r \in \{ p, \ldots, \left\lceil \frac{n}{2} \right\rceil +1 \}$ in polynomial time, and, similarly, if there are integers $p$ and $q$ with $p \le q$ such that one can find spanning trees of a connected graph $G'$ with $p$ and $q$ leaves, respectively, then one can find a spanning tree of $G'$ with $r$ leaves for every integer $r \in \{ p, \ldots, q \}$ in polynomial time. 
We now sketch a proof and leave the details to the reader.

As mentioned in the introduction, one can infer from Theorem A~\cite{CRS12} that every connected and locally connected $n$-vertex graph $G$ has a spanning tree with at least $\frac{n}{2} + 1$ leaves; it follows from their proof, in particular their Lemmas 2.3 and 2.4, that this can be done in polynomial time (one first builds a spanning weak 2-tree in $G$ and then a HIST based on the weak 2-tree). The above statements now follow from our theorem and its proof as follows. In Claim~1, we need to perform at most $n-1$ edge exchanges, and in Claim 2, we need to do at most $\frac{n}{2}$ edge exchanges to find a spanning tree $S_i$. So we need to perform at most $\left\lfloor \frac{n}{2} \right\rfloor \cdot (n-1)$ edge exchanges in total.

For planar 4-connected graphs, we remark that a hamiltonian cycle in such graphs can be found in linear time as proven by Chiba and Nishizeki~\cite{CN89}. So in any planar 4-connected graph $G$, for any integer $k$ with $2 \le k \le \left\lceil \frac{|V(G)|}{2} \right\rceil + 1$ we can find a spanning tree of $G$ with exactly $k$ leaves in polynomial time.

We have already discussed a lower bound for $f(n)$. We now present an upper bound. 

\medskip

\noindent {\bf Proposition.} {\it There exist infinitely many $n$ such that there is a plane $4$-connected $n$-vertex triangulation containing a spanning tree with $\frac{2n}{3}$ leaves, but no spanning tree with more than $\frac{2n}{3}$ leaves.}

\begin{proof}
The following construction is illustrated in Fig.~\ref{fig:1}. Let $C_4$ be the cycle on 4 vertices and $P_k$ the path on $k$ vertices, where $k \ge 6$ shall be a multiple of 3. Consider the cartesian product of $C_4$ and $P_k$. 
This graph contains $k$ pairwise disjoint 4-cycles $C^1, \ldots, C^k$, where $C^i = v^i_0v^i_1v^i_2v^i_3v^i_0$ for $i \in \{ 1, \ldots, k \}$, such that $v^i_j$ is adjacent to $v^{i+1}_j$ for every pair of $i \in \{ 1, \ldots, k-1 \}$ and $j \in \{ 0, 1, 2, 3 \}$. 
We add, for every $i \in \{ 1, \ldots, k - 1 \}$, the edge $v^i_jv^{i+1}_{j+1}$ to this graph, for every $j \in \{ 0, 1, 2, 3 \}$, indices taken modulo~4, and also add the edges $v^1_1v^1_3$ and $v^k_1v^k_3$. 
Denote the resulting graph, which is a plane 4-connected triangulation of order $4k =: n$, by $G_k$.

Consider, for $i \in \{ 1, \ldots, k - 2\}$, the subgraph of $G_k$ induced by $V(C^i) \cup V(C^{i+1}) \cup V(C^{i+2})$ and denote it by $H_i$. We see $C^i, C^{i+1}, C^{i+2}$ as subgraphs of $H_i$. For a tree $T$, we call {\it stem vertices} the vertices in $T$ that are not leaves of $T$. Let $T$ be a spanning tree of $G_k$ and $S$ the subtree of $T$ formed by $T$'s stem vertices. Since $S$ is connected and also by the structure of $G_k$, $S \cap H_i$ must contain a path $P$ with at least one vertex in $C^i$ and at least one vertex in $C^{i+2}$ for all $i \in \{ 2, \ldots, k - 3\}$. Thus $P$ has at least three vertices. Assume it has exactly three and put $P = uvw$. There is a vertex $v'$ in $C^{i+1}$ at distance 2 from $v$ in $G_k$. 
It can be easily seen that no matter the position of $u$ and $w$ in $H_i$, they will never be adjacent to $v'$. We call this argument ($\dagger$). Since $T$ is a spanning tree of $G_k$, by ($\dagger$) we must have $|V(S) \cap V(H_i)| \ge 4$ for all $i \in \{ 2, \ldots, k - 3\}$.

If $V(S) \cap V(C^1) \ne \emptyset$, by ($\dagger$) we have $|V(S) \cap V(H_1)| \ge 4$. Assume henceforth $V(S) \cap V(C^1) = \emptyset$. Since the vertices $v^1_1$ and $v^1_3$ lie in $V(T) \setminus V(S)$, each of them must be adjacent to a vertex in $S$. By the structure of $G_k$, this cannot be the same vertex, so $|V(S) \cap V(C^2)| \ge 2$. If $|V(S) \cap V(C^2)| = 2$ and the two vertices in $V(S) \cap V(C^2)$ are adjacent, again by the structure of $G_k$ it is clear that not all vertices of $C^1$ are adjacent to some vertex in $S$, a contradiction. So in this case necessarily $|V(S) \cap V(C^2)| \ge 3$, and as $k\ge 6$ and by the connectedness of $S$ we obtain $|V(S) \cap V(H_1)| \ge 4$. 
If the two vertices are non-adjacent, as $S$ is connected and by the structure of $G_k$, we have $|V(S) \cap V(C^3)| \ge 2$. So once more $|V(S) \cap V(H_1)| \ge 4$. Analogously one proves that $|V(S) \cap V(H_{k-2})| \ge 4$. We have proven that $|V(S) \cap V(H_i)| \ge 4$ for all $i \in \{ 1, \ldots, k - 2\}$.

\begin{figure}
  \begin{center}
    \includegraphics[height=55mm]{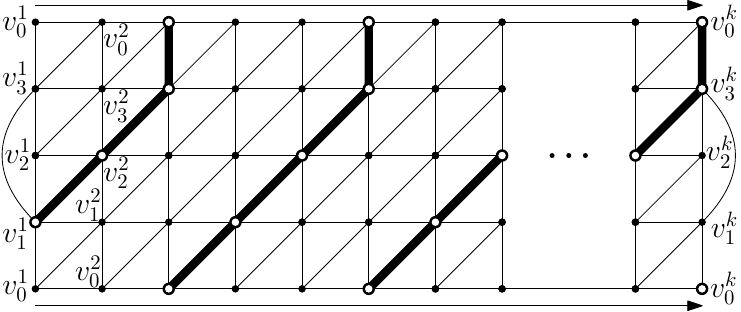}
    \caption{A plane 4-connected triangulation $G_k$ obtained by identifying the top and bottom paths using the given orientation. The stem vertices of a spanning tree are emphasised.}
    \label{fig:1}
  \end{center}
\end{figure}

Since $|V(S) \cap V(H_i)| \ge 4$ for all $i \in \{ 3j+1 \mid j\in\{0, \ldots, \frac{k}{3}-1\}\}$ and $V(G_k)$ is partitioned into $V(H_1) \cup V(H_4) \cup \cdots \cup V(H_{k-2})$, 
we can conclude that every spanning tree in $G_k$ has at least $\frac{n}{3}$ stem vertices, and thus has at most $\frac{2n}{3}$ leaves. In Fig.~\ref{fig:1}, a path $Q$ on $\frac{n}{3}$ vertices in $G_k$ is shown by bold lines. As $k$ is a multiple of 3, every vertex in $G_k$ either lies in $Q$ or is adjacent to a vertex in $Q$. We therefore obtain from $Q$ a spanning tree of $G_k$ with $\frac{2n}{3}$ leaves. 
\end{proof}

\noindent \textbf{Remark 1.} The upper bound $\frac{2n}{3}$ (perhaps plus some constant) can be seen as stemming from $6$-regularity. For a $6$-regular triangulation $G$ (of the torus or the Klein bottle), every subtree of $G$ with $k$ stems has at most $2k+4$ leaves as shown below. 
Thus, for any spanning tree of $G$, the number of leaves is at most $\frac{(2k+4)n}{3k+4}$.

Let $T$ be a subtree of $G$ with $k$ stem vertices and $S$ be the set of all its stem vertices. Take an arbitrary sequence of connected subgraphs $S=S_k \supseteq S_{k-1} \supseteq \cdots \supseteq S_1 =: \{ v_1 \}$ where $S_i$ consists of exactly $i$ vertices. 
We show by induction that for every $1\le i \le k$, every subtree of $G$ in which the set of stem vertices is $S_i$ has at most $2i+4$ leaves.

\newpage

For $i=1$, it is trivial. Suppose that the statement is true for $i-1$. Fix a subtree $T_{i-1}$ with the stem vertices $S_{i-1}$ such that the number of leaves is maximum. 
Let $T_i$ be a subtree of $G$ with the stem vertices $S_i$. Now $v_i:= V(S_i)\backslash V(S_{i-1})$ should be a leaf in $T_{i-1}$ since $v_i$ has a neighbour $v_j$ in $S_{i-1}$ for some $j<i$. Note that $v_i$ and $v_j$ have two common neighbours, say $x$ and $y$, both of which are in $T_{i-1}$ (as leaves or stem vertices). Then, $V(T_i)\backslash V(T_{i-1}) \subseteq N_G(v_i)\backslash \{v_j, x, y\}$ and $|V(T_i)\backslash V(T_{i-1})| \le 3$. 
Thus, the number of leaves in $T_i$ is at most $(2(i-1)+4)+(3-1) = 2i+4$. 

\bigskip

Our lower and upper bounds on $f(n)$ lie far apart, and we do not know how to close this gap. We only have limited evidence, but we suspect that the upper bound $\frac{2n}{3}$ is tight, which was also raised as a question by Bradshaw et al.~\cite[Question 4.3]{BMNS22} recently. 

\bigskip

\noindent \textbf{Remark 2.} Very recently, Bose et al.~\cite{BDHM23+} announced that they showed the lower bound $\frac{11n}{21}$ for plane (not necessarily 4-connected) $n$-vertex triangulations. This yields an improvement of the lower bound for $f(n)$. Still, the gap between it and the upper bound given here remains large.

\bigskip


\label{sec:biblio}

\end{document}